\documentclass[twoside]{article}%
\usepackage{amssymb}
\usepackage{amsfonts}
\usepackage{amsmath}
\usepackage{graphicx}%
\providecommand{\U}[1]{\protect\rule{.1in}{.1in}}
\newtheorem{theorem}{Theorem}

\newtheorem{proposition}[theorem]{Proposition}
\newtheorem{remark}[theorem]{Remark}

\newenvironment{proof}[1][Proof]{\noindent\textbf{#1.} }{\ \rule{0.5em}{0.5em}}
\begin{document}

\title{\textbf{Blowup of }$C^{2}$\textbf{ Solutions for the Euler Equations and
Euler-Poisson Equations in }$R^{N}$}
\author{Y\textsc{uen} M\textsc{anwai\thanks{E-mail address: nevetsyuen@hotmail.com }}\\\textit{Department of Applied Mathematics,}\\\textit{The Hong Kong Polytechnic University,}\\\textit{Hung Hom, Kowloon, Hong Kong}}
\date{Revised 12-Sept-2009}
\maketitle

\begin{abstract}
In this paper, we use integration method to show that there is no existence of
global $C^{2}$ solution with compact support, to the pressureless
Euler-Poisson equations with attractive forces in $R^{N}$. And the similar
result can be shown, provided that\ the uniformly bounded functional:%
\begin{equation}
\int_{\Omega(t)}K\gamma(\gamma-1)\rho^{\gamma-2}(\nabla\rho)^{2}%
dx+\int_{\Omega(t)}K\gamma\rho^{\gamma-1}\Delta\rho dx+\epsilon\geq
-\delta\alpha(N)M,
\end{equation}
where $M$ is the mass of the solutions and $\left\vert \Omega\right\vert $ is
the fixed volume of $\Omega(t)$.

On the other hand, our differentiation method provides a simpler proof to show
the blowup result in "D. H. Chae and E. Tadmor, \textit{On the Finite Time
Blow-up of the Euler-Poisson Equations in }$R^{N}$, Commun. Math. Sci.
\textbf{6} (2008), no. 3, 785--789.".


\end{abstract}

\section{Introduction}

The Euler $(\delta=0)$/ Euler-Poisson $(\delta=\pm1)$ equations can be written
in the following form:
\begin{equation}
\left\{
\begin{array}
[c]{rl}%
{\normalsize \rho}_{t}{\normalsize +\nabla\cdot(\rho u)} & {\normalsize =}%
{\normalsize 0,}\\
{\normalsize (\rho u)}_{t}{\normalsize +\nabla\cdot(\rho u\otimes u)+\nabla
}P & {\normalsize =}{\normalsize \delta\rho\nabla\Phi,}\\
{\normalsize \Delta\Phi(t,x)} & {\normalsize =\alpha(N)}{\normalsize \rho,}%
\end{array}
\right.  \label{Euler-Poisson}%
\end{equation}
where $\alpha(N)$ is a constant related to the unit ball in $R^{N}$:
$\alpha(1)=1,$ $\alpha(2)=2\pi$ and $\alpha(3)=4\pi$. And as usual, $\rho
=\rho(t,x)$ and $u=u(t,x)\in\mathbf{R}^{N}$ are the density and the velocity
respectively. $P=P(\rho)$\ is the pressure function. In the above systems, the
self-gravitational potential field $\Phi=\Phi(t,x)$\ is determined by the
density $\rho$ itself, through the Poisson equation (\ref{Euler-Poisson}%
)$_{3}$. For $N=3$, the equations (\ref{Euler-Poisson}) are the classical
(non-relativistic) descriptions of a galaxy, in astrophysics. See \cite{BT}
and \cite{C}, for details about the systems. The $\gamma$-law can be applied
on the pressure term $P(\rho)$, i.e.%
\begin{equation}
{\normalsize P}\left(  \rho\right)  {\normalsize =K\rho}^{\gamma},
\label{gamma}%
\end{equation}
which is a common hypothesis. If the parameter\ ${\normalsize K>0}$, we call
the system with pressure; if ${\normalsize K=0}$, we call it pressureless. The
constant $\gamma=c_{P}/c_{v}\geq1$, where $c_{P}$, $c_{v}$\ are the specific
heats per unit mass under constant pressure and constant volume respectively,
is the ratio of the specific heats, that is, the adiabatic exponent in the
equation (\ref{gamma}). In particular, the fluid is called isothermal if
$\gamma=1$. If the parameter\ ${\normalsize K>0}$, we call the system with
pressure; if ${\normalsize K=0}$, we call it pressureless.\newline On the
other hand, the Poisson equation (\ref{Euler-Poisson})$_{3}$ can be solved as%
\begin{equation}
{\normalsize \Phi(t,x)=}\int_{R^{N}}G(x-y)\rho(t,y){\normalsize dy,}%
\end{equation}
where $G$ is the Green's function for the Poisson equation in the
$N$-dimensional spaces defined by
\begin{equation}
G(x)\doteq\left\{
\begin{array}
[c]{ll}%
|x|, & N=1;\\
\log|x|, & N=2;\\
\frac{-1}{|x|^{N-2}}, & N\geq3.
\end{array}
\right.
\end{equation}
When $\delta=1$, the system is the compressible Euler equations with repulsive
forces. The equation (\ref{Euler-Poisson})$_{3}$ is the Poisson equation
through which the potential with repulsive forces is determined by the density
distribution of the electrons. In this case, the equations can be viewed as a
semiconductor model. See \cite{Cse}, \cite{Lions} for details about the
system. When $\delta=-1$, the system can model the self-gravitating fluid,
such as gaseous stars. Besides, the evolution of the cosmology can be modelled
by the dust distribution without pressure term. That describes the stellar
systems of collisionless and gravitational $n$-body systems \cite{FT}. And the
pressureless Euler-Poisson equations can be derived from the
Vlasov-Poisson-Boltzmann model with the zero mean free path \cite{G1}.\newline
Usually the Euler-Poisson equations can be rewritten in the scalar form:%
\begin{equation}
\left\{
\begin{array}
[c]{rl}%
\frac{\partial\rho}{\partial t}+\underset{k=1}{\overset{N}{\Sigma}}u_{k}%
\frac{\partial\rho}{\partial x_{k}}+\rho\underset{k=1}{\overset{N}{\Sigma}%
}\frac{\partial u_{k}}{\partial x_{k}} & {\normalsize =}{\normalsize 0,}\\
\rho\left(  \frac{\partial u_{i}}{\partial t}+\underset{k=1}{\overset
{N}{\Sigma}}u_{k}\frac{\partial u_{i}}{\partial x_{k}}\right)  +\frac{\partial
P}{\partial x_{i}} & {\normalsize =\delta\rho\frac{\partial\Phi}{\partial
x_{i}}}\text{, for }i=1,2,...N.
\end{array}
\right.  \label{gamma=1}%
\end{equation}
For the construction of the analytical solutions for the system with pressure,
the interested readers may see \cite{C}, \cite{GW}, \cite{M}, \cite{DXY},
\cite{CT}, \cite{Y}, \cite{Y1}, \cite{Y3} and \cite{Y4}. The results for local
existence theories can be founded in \cite{M}, \cite{B} and \cite{G}. The
analysis of stabilities for the systems may be referred in \cite{SI},
\cite{A}, \cite{H}, \cite{E}, \cite{MUK}, \cite{MP}, \cite{P}, \cite{DLY},
\cite{DXY}, \cite{Y1}, \cite{J}, \cite{Y2}, \cite{CT} and \cite{CH}.

Recently, Chae and Tadmor \cite{CT} showed the finite time blowup, for the
pressureless Euler-Poisson equations with attractive forces, under the initial
condition,%
\begin{equation}
S:=\{\left.  a\in R^{N}\right\vert \text{ }\rho_{0}(a)>0,\text{ }%
\Omega(a)=0,\text{ }\nabla\cdot u(0,x(0))<0\}\neq\phi.\label{chea}%
\end{equation}
On the other hand, in \cite{Y2}, we have the blowup results if the solutions
with compact support under the condition,
\begin{equation}
2\int_{\Omega(t)}(\rho\left\vert u\right\vert ^{2}+2P)dx<M^{2}-\epsilon,
\end{equation}
where $M$ is the mass of the solution.

In this article, the alternative approaches are adopted to show that there is
no global existence of $C^{2}$ solutions for the system, (\ref{gamma=1}), with
compact support without the condition (\ref{chea}):

\begin{theorem}
\label{thm:1 copy(1)}For the pressureless Euler-Poisson equations with
attractive forces $(\delta=-1)$, (\ref{gamma=1}), there do not exist global
$C^{2}$ solutions $(\rho,u)$ with compact support. For the systems with
pressure, for $\gamma>1$, the above result is also true provided that the
uniformly bounded functional:%
\begin{equation}
\int_{\Omega(t)}K\gamma\left[  (\gamma-1)\rho^{\gamma-2}(\nabla\rho
)^{2}dx+\rho^{\gamma-1}\Delta\rho\right]  dx+\epsilon\geq-\delta
\alpha(N)M,\label{condi1}%
\end{equation}
where $\epsilon$ is an arbitrary small positive constant, $M$ is the mass of
the solution and $\left\vert \Omega\right\vert $ is the fixed volume of
$\Omega(t)$.
\end{theorem}

\section{Integration Method}

In this section, we present the proof of Theorem \ref{thm:1 copy(1)}.

\begin{proof}
First, we show that the $\rho(t,x(t;x))$ preserves its positive nature as the
mass equation (\ref{gamma=1})$_{1}$ can be converted to be%
\begin{equation}
\frac{D\rho}{Dt}+\rho\nabla\cdot u=0, \label{io1}%
\end{equation}
with the material derivative,%
\begin{equation}
\frac{D}{Dt}=\frac{\partial}{\partial t}+\left(  u\cdot\nabla\right)  .
\end{equation}
Integrate the equation (\ref{io1})$:$%
\begin{equation}
\rho(t,x)=\rho_{0}(x_{0}(0,x_{0}))\exp\left(  -\int_{0}^{t}\nabla\cdot
u(t,x(t;0,x_{0}))dt\right)  >0,
\end{equation}
for $\rho_{0}(x_{0}(0,x_{0}))>0,$along the characteristic curve.

For the pressureless Euler-Poisson equations with attractive forces
$(\delta=-1)$, we divide $\rho$ to the momentum equation (\ref{gamma=1})$_{2}$
to have%
\begin{equation}
u_{t}+u\cdot\nabla u=-\nabla\Phi,
\end{equation}
And we take divergence to the above equation:
\begin{equation}
\nabla(u_{t}+u\cdot\nabla u)=-\Delta\Phi.
\end{equation}
By taking integration to the above equation, we have,%
\begin{equation}
\int_{\Omega(t)}\nabla(u_{t}+u\cdot\nabla u)dx=-\int_{\Omega(t)}\alpha(N)\rho
dx,
\end{equation}%
\begin{equation}
\int_{\Omega(t)}\nabla(u_{t}+u\cdot\nabla u)dx=-\alpha(N)M,
\end{equation}%
\begin{equation}
\int_{\Omega(t)}\left[  \underset{i=1}{\overset{N}{\Sigma}}u_{itx_{i}%
}+\underset{i=1}{\overset{N}{\Sigma}}u_{i}\left(  \underset{k=1}{\overset
{N}{\Sigma}}u_{kx_{i}x_{k}}\right)  +\underset{i=1}{\overset{N}{\Sigma}%
}u_{ix_{i}}^{2}\right]  dx=-\alpha(N)M.
\end{equation}
It becomes alone the characteristic curve:%
\begin{equation}
\int_{R^{N}}\frac{D}{Dt}\nabla\cdot u(t,x(t;x))dx+\int_{\Omega(t)}\left[
\nabla\cdot u(t,x(t;x))\right]  ^{2}dx\leq-\alpha(N)M.
\end{equation}
By denoting
\begin{equation}
H:=H(t,x)=\int_{\Omega(t)}\underset{i=1}{\overset{N}{\Sigma}}u_{ix_{i}%
}(t,x(t;x))dx=\int_{\Omega(t)}\nabla\cdot u(t,x(t;x))dx,
\end{equation}
and with the Cauchy-Schwarz inequality,
\begin{align}
\left\vert \int_{\Omega(t)}\nabla\cdot u(t,x(t;x))\cdot1dx\right\vert  &
\leq\left(  \int_{\Omega(t)}\left[  \nabla\cdot u(t,x(t;x))\right]
^{2}dx\right)  ^{1/2}\left(  \int_{\Omega(t)}1dx\right)  ^{1/2},\\
\frac{\left\vert \int_{\Omega(t)}\nabla\cdot u(t,x(t;x))dx\right\vert
}{\left\vert \Omega\right\vert ^{1/2}} &  \leq\left(  \int_{\Omega(t)}\left[
\nabla\cdot u(t,x(t;x))\right]  ^{2}dx\right)  ^{1/2},\\
\frac{H(t)^{2}}{\left\vert \Omega\right\vert } &  \leq\int_{\Omega(t)}\left[
\nabla\cdot u(t,x(t;x))\right]  ^{2}dx,
\end{align}
where $\left\vert \Omega\right\vert $ is the fixed volume of $\Omega
(t)$,\newline we have%
\begin{equation}
\frac{DH}{Dt}+\frac{H^{2}}{\left\vert \Omega\right\vert }\leq-\alpha(N)M,
\end{equation}%
\begin{equation}
H(t)\leq-\sqrt{\left\vert \Omega\right\vert \alpha(N)M}\tan\left(  \sqrt
{\frac{\alpha(N)M}{\left\vert \Omega\right\vert }}t-\tan^{-1}\left(
\sqrt{\frac{1}{\left\vert \Omega\right\vert \alpha(N)M}}H_{0}\right)  \right)
,
\end{equation}%
\begin{equation}
H(T)\leq-\infty,
\end{equation}
with the finite time $T>0$, such that%
\begin{equation}
\sqrt{\frac{\alpha(N)M}{\left\vert \Omega\right\vert }}T=\tan^{-1}\left(
\sqrt{\frac{1}{\left\vert \Omega\right\vert \alpha(N)M}}H_{0}\right)
+\frac{\pi}{2},
\end{equation}%
\begin{equation}
T=\sqrt{\frac{\left\vert \Omega\right\vert }{\alpha(N)M}}\left[  \tan
^{-1}\left(  \sqrt{\frac{1}{\left\vert \Omega\right\vert \alpha(N)M}}%
H_{0}\right)  +\frac{\pi}{2}\right]  .
\end{equation}
Therefore, for any $C^{2}$ solutions with compact support, they blow up before
$T=\sqrt{\frac{\left\vert \Omega\right\vert }{\alpha(N)M}}\pi$.

On the other hand, for the systems with pressure $(\gamma>1)$, we divide
$\rho$ to the momentum equation (\ref{gamma=1})$_{2}$ to have,%
\begin{equation}%
\begin{array}
[c]{rl}%
\left(  \frac{\partial u_{i}}{\partial t}+\underset{k=1}{\overset{N}{\Sigma}%
}u_{k}\frac{\partial u_{i}}{\partial x_{k}}\right)  +K\gamma\rho^{\gamma
-1}\nabla\rho & {\normalsize =\delta\frac{\partial\Phi}{\partial x_{i}}%
}\text{.}%
\end{array}
\end{equation}
Take differentiation to the momentum equation for $\gamma>1$:%
\begin{align}
\nabla(u_{t}+u\cdot\nabla u)+K\gamma(\gamma-1)\rho^{\gamma-2}(\nabla\rho
)^{2}+K\gamma\rho^{\gamma-1}\Delta\rho &  =-\delta\Delta\Phi,\\
\nabla(u_{t}+u\cdot\nabla u)+K\gamma(\gamma-1)\rho^{\gamma-2}(\nabla\rho
)^{2}+K\gamma\rho^{\gamma-1}\Delta\rho &  =-\delta\alpha(N)\rho,\nonumber
\end{align}
with the Poisson equation (\ref{Euler-Poisson})$_{3}$. By taking integration
to the above equation, we have,%
\begin{equation}
\int_{\Omega(t)}\nabla(u_{t}+u\cdot\nabla u)dx+\int_{\Omega(t)}\left[
K\gamma(\gamma-1)\rho^{\gamma-2}(\nabla\rho)^{2}+K\gamma\rho^{\gamma-1}%
\Delta\rho\right]  dx=-\int_{\Omega(t)}\delta\alpha(N)\rho dx,
\end{equation}%
\begin{equation}
\int_{\Omega(t)}\nabla(u_{t}+u\cdot\nabla u)dx+\int_{\Omega(t)}K\gamma
(\gamma-1)\rho^{\gamma-2}(\nabla\rho)^{2}dx+\int_{\Omega(t)}K\gamma
\rho^{\gamma-1}\Delta\rho dx=-\delta\alpha(N)M,
\end{equation}%
\begin{equation}
\int_{\Omega(t)}\left[  \underset{i=1}{\overset{N}{\Sigma}}u_{itx_{i}%
}+\underset{i=1}{\overset{N}{\Sigma}}u_{i}\left(  \underset{k=1}{\overset
{N}{\Sigma}}u_{kx_{i}x_{k}}\right)  +\underset{i=1}{\overset{N}{\Sigma}%
}u_{ix_{i}}^{2}\right]  dx\leq-\epsilon,
\end{equation}
with the required condition (\ref{condi1}).\newline Then for any $H(0)$, we
have:%
\begin{align}
\frac{DH(t)}{Dt}+\frac{H(t)^{2}}{\left\vert \Omega\right\vert } &
\leq-\epsilon,\\
H(T) &  \leq-\infty,
\end{align}
with a finite time $T$.\newline Therefore, for any $C^{2}$ solutions with
compact support, they blow up on or before a finite time $T$.

This completes the proof.
\end{proof}

\section{Differentiation Method}

Alternatively, we may take differentiation to the momentum equations
(\ref{gamma=1})$_{2}$ to have the same result of in \cite{CT}:

\begin{proposition}
Suppose $(\rho,u)$ are $C^{2}$ solutions for the pressureless $(\delta=0)$
Euler-Poisson equations in $R^{N}$, (\ref{gamma=1}), the solutions blow up
before $T=-1/H_{0}$, with the initial velocity at some non-vacuum point:%
\begin{equation}
H_{0}=\nabla\cdot u(0,x_{0})<0\text{.} \label{add2}%
\end{equation}

\end{proposition}

\begin{proof}
For $\rho_{0}(x_{0}(0,x_{0}))>0$, along the characteristic curve, we have:%
\begin{align}
u_{t}+u\cdot\nabla u &  =-\nabla\Phi,\\
\nabla(u_{t}+u\cdot\nabla u) &  =-\Delta\Phi,\\
\underset{i=1}{\overset{N}{\Sigma}}u_{itx_{i}}+\underset{i=1}{\overset
{N}{\Sigma}}u_{i}\left(  \underset{k=1}{\overset{N}{\Sigma}}u_{kx_{k}x_{i}%
}\right)  +\frac{1}{N}\left(  \underset{i=1}{\overset{N}{\Sigma}}u_{ix_{i}%
}\right)  ^{2}+\alpha(N)\rho &  =0,
\end{align}
with the required condition (\ref{add2}). By defining $H=H(t):=\underset
{i=1}{\overset{N}{\Sigma}}u_{ix_{i}}(t,x(t))$, we have%
\begin{equation}
\frac{DH}{Dt}+\frac{1}{N}H^{2}+\alpha(N)\rho=0,\label{we1}%
\end{equation}
with the Poisson equation (\ref{Euler-Poisson})$_{3}$. It becomes
\begin{align}
\frac{DH}{Dt} &  \leq-\frac{1}{N}H^{2},\\
H(t) &  \leq\frac{N\underset{i=1}{\overset{N}{\Sigma}}u_{ix_{i}}(0;(0,x_{0}%
))}{N+\underset{i=1}{\overset{N}{\Sigma}}u_{ix_{i}}(0;(0,x_{0}))t}%
\rightarrow-\infty,
\end{align}
as $t\rightarrow-N/\left[  \underset{i=1}{\overset{N}{\Sigma}}u_{ix_{i}%
}(0;(0,x_{0}))\right]  :=-N/\left[  \underset{i=1}{\overset{N}{\Sigma}%
}u_{ix_{i}}(0,x_{0})\right]  $ and $\underset{i=1}{\overset{N}{\Sigma}%
}u_{ix_{i}}(0;(0,x_{0}))<0,$ for some point $x_{0}$. Therefore, we show that
there does not exist global $C^{2}$ solutions, with the initial velocity at
some point,%
\begin{equation}
H_{0}=\nabla\cdot u(0,x_{0})<0.
\end{equation}

The proof is completed.
\end{proof}

\begin{remark}
The above method provides a simpler way to show the show the same inequality
in \cite{CT},%
\begin{equation}
\frac{D\operatorname{div}u}{Dt}\leq-\frac{1}{N}(\operatorname{div}u)^{2},
\end{equation}
without the analysis of spectral dynamics.
\end{remark}

\section{Discussion}

Makino, Ukai, Kawashima initially defined the "Tame" solutions (the solutions
to the pressureless Euler equations) \cite{MUK}. Then Makino and Perthame
considered the Tame solutions for the system with gravitational forces
\cite{MP}. After that Perthame studied the $3$-dimensional pressureless system
with repulsive forces \cite{P}. All of the above results rely on the solutions
with radial symmetry:
\begin{equation}%
\begin{array}
[c]{rl}%
u_{t}+uu_{r} & {\normalsize =}\frac{\alpha(N)\delta}{r^{N-1}}\int_{0}^{r}%
\rho(t,s)s^{N-1}ds,
\end{array}
\label{ss1}%
\end{equation}
where $r:=\left(  \sum_{i=1}^{N}x_{i}^{2}\right)  ^{1/2}$ is the radius. And
the Emden ordinary differential equations were deduced on the boundary point
of the solutions with compact support:%
\begin{equation}
\frac{D^{2}R}{Dt^{2}}=\frac{\delta M}{R^{N-1}},\text{ }R(0,R_{0})=R_{0}%
\geq0,\text{ }\dot{R}(0,R_{0})=0,
\end{equation}
where $\frac{dR}{dt}:=u$ and $M$ is the mass of the solutions, along the
characteristic curve. They showed the blowup results for the $C^{1}$ solutions
of the system (\ref{ss1}). And recently, Chae and Tadmor \cite{CT} obtain the
blowup result, which does not require the solutions in radial symmetry.
However, all the above results concern about the pressureless cases with the
external forces only. This article has shed new light on the situations with
the pressure term. In particular, it answers some cases for the Euler
equations or the Euler Poisson equations. However, the condition
(\ref{condi1}) in our theorem is too restricted. An refinement for it is
expected in the future work.


\begin{thebibliography}{99}                                                                                               %


\bibitem {A}S. Alinhac, \textit{Blowup for Nonlinear Hyperbolic Equations}.
Progress in Nonlinear Differential Equations and their Applications, 17.
Birkh\"{a}ser Boston, Inc., Boston, MA, 1995.

\bibitem {B}M. Bezard, \textit{Existence locale de solutions pour les
equations d'Euler--Poisson (Local Existence of Solutions for Euler--Poisson
Equations)}, Japan J. Indust. Appl. Math. \textbf{10} (3) (1993) 431--450 (in French).

\bibitem {BT}J. Binney and S. Tremaine, \textit{Galactic Dynamics}, Princeton
Univ. Press, 1994.

\bibitem {Cse}C. F. Chen, \textit{Introduction to Plasma Physics and
Controlled Fusion}, Plenum, New York (1984)

\bibitem {CT}D. H. Chae and E. Tadmor, \textit{On the Finite Time Blow-up of
the Euler-Poisson Equations in }$R^{N}$, Commun. Math. Sci. \textbf{6} (2008),
no. 3, 785--789.

\bibitem {CH}D. H. Chae and S. Y. Ha, \textit{On the Formation of Shochs to
the Compressible Euler Equations}, Article in Press, Commun. Math. Sci.

\bibitem {C}S. Chandrasekhar, \textit{An Introduction to the Study of Stellar
Structure}, Univ. of Chicago Press, 1939.

\bibitem {DLY}Y.B. Deng, T.P. Liu, T. Yang and Z.A. Yao, \textit{Solutions of
Euler-Poisson Equations for Gaseous Stars}, Arch. Ration. Mech. Anal.
\textbf{164} (2002) 261--285.

\bibitem {DXY}Y.B. Deng, J.L. Xiang, T. Yang,\textit{ Blowup Phenomena of
Solutions to Euler-Poisson Equations}, J. Math. Anal. Appl. \textbf{286} (1)
(2003) 295--306.

\bibitem {E}S. Engelberg, \textit{Formation of Singularities in the Euler and
Euler-Poisson Equations}, Phys. D, 98(1), 67-74, 1996.

\bibitem {FT}H. H. Fliche and R. Triay, \textit{Euler-Poisson-Newton Approach
in Cosmology}, Cosmology and Gravitation, 346--360, AIP Conf. Proc., 910,
Amer. Inst. Phys., Melville, NY, 2007.

\bibitem {G}P. Gamblin, \textit{Solution reguliere a temps petit pour
l'equation d'Euler--Poisson (Small-time Regular Solution for the
Euler--Poisson Equation)}, Comm. Partial Differential Equations \textbf{18}
(5--6) (1993) 731--745 (in French). Comm. Partial Differential Equations
\textbf{18} (5--6) (1993) 731--745 (in French).

\bibitem {G1}R. T. Glassey, \textit{The Cauchy Problem in Kinetic Theory},
Society for Industrial and Applied Mathematics (SIAM), Philadelphia, PA, 1996.

\bibitem {GW}P. Goldreich and S. Weber, \textit{Homologously Collapsing
Stellar Cores}, Astrophys, J. \textbf{238} (1980), 991-997 .

\bibitem {H}L. H\"{o}rmander, \textit{Lectures on Non-linear Hyperbolic
Differential Equations}, Mathematics and Applications, Vol 26, Springer, 1997.

\bibitem {J}J. Jang, \textit{Nonlinear Instability in Gravitational
Euler-Poisson Systems for} $\gamma=6/5$, Arch. Ration. Mech. Anal.
\textbf{188} (2)(2008), 265--307.

\bibitem {Li}T.H. Li, \textit{Some Special Solutions of the Multidimensional
Euler Equations in }$R^{N}$, Comm. Pure Appl. Anal\textit{.} \textbf{4}
(4)(2005), 757-762.

\bibitem {Lions}P.L. Lions, \textit{Mathematical Topics in Fluid Mechanics}.
Vols. 1, 2, 1998, Oxford: Clarendon Press, 1998.

\bibitem {M}T. Makino, \textit{On a Local Existence Theorem for the Evolution
Equation of Gaseous Stars}. Patterns and waves, 459--479, Stud. Math. Appl.,
\textbf{18}, North-Holland, Amsterdam, 1986.

\bibitem {MUK}T. Makino, S. Ukai and S. Kawashima, \textit{On Compactly
Supported Solutions of the Compressible Euler Equation}, Recent Topics in
Nonlinear PDE, III (Tokyo, 1986), 173--183, North-Holland Math. Stud., 148,
North-Holland, Amsterdam, 1987.

\bibitem {MP}T. Makino and B. Perthame, \textit{Sur les Solutions a symmetric
spherique de lequation d'Euler-poisson Pour levolution d'etoiles
gazeuses,(French) [On Radially Symmetric Solutions of the Euler-Poisson
Equation for the Evolution of Gaseous Stars], }Japan J. Appl. Math. \textbf{7}
(1990), 165-170.

\bibitem {P}B. Perthame, \textit{Nonexistence of Global Solutions to
Euler-Poisson Equations for Repulsive Forces}, Japan J. Appl. Math. \textbf{7}
(1990), 363--367.

\bibitem {SI}T.C. Sideris, \textit{Formation of Singularities in
Three-dimensional Compressible Fluids}, Comm. Math. Phys. \textbf{101} (1985),
No. 4, 475-485.

\bibitem {Y}M.W. Yuen, \textit{Blowup Solutions for a Class of Fluid Dynamical
Equations in }$R^{N}$, J. Math. Anal. Appl. \textbf{329} (2)(2007), 1064-1079.

\bibitem {Y1}M.W. Yuen, \textit{Analytical Blowup Solutions to the }%
$2$\textit{-dimensional Isothermal Euler-Poisson Equations of Gaseous Stars},
J. Math. Anal. Appl. \textbf{341 (}1\textbf{)(}2008\textbf{), }445-456.

\bibitem {Y2}M.W. Yuen, \textit{Stabilities for Euler-Poisson Equations in
Some Special Dimensions}, J. Math. Anal. Appl. \textbf{344} \textbf{(}%
1\textbf{)}(2008), 145--156.

\bibitem {Y3}M.W. Yuen, \textit{Analytical Blowup Solutions to the }%
$2$\textit{-dimensional Isothermal Euler-Poisson Equations of Gaseous Stars
}$II$, preprint, arXiv:0906.0176v2.

\bibitem {Y4}M.W. Yuen, \textit{Analytical Blowup Solutions to the Isothermal
Euler-Poisson Equations of Gaseous Stars in }$R^{N}$, preprint, arXiv:0906.0178v1.
\end{thebibliography}
\end{document}